# REAL AND COMPLEX ORDER INTEGRALS AND DERIVATIVES OPERATORS OVER THE SET OF CAUSAL FUNCTIONS

**Raoelina Andriambololona**
[1] Theoretical Physics Department, Institut National des Sciences et Techniques Nucléaires (INSTN-Madagascar)
raoelinasp@yahoo.fr, jacquelineraoelina@hotmail.com

*Abstract*— The fractional integrals and fractional derivatives problem is tackled by using the operator approach. The definition domain $E$ of operators is the set of causal functions. Many properties of fractional integrals are given. Fractional derivatives definition is derived from fractional integrals one. Then an unified definitions of fractional integrals and derivatives operator is obtained according to the sign of the real part of the order $s$. The study utilizes many properties of the Euler's gamma aand beta functions and their extensions in $\mathbb{R}$ and $\mathbb{C}$. Comparison with the definitions given by other authors (Liouville, Riemann, Liouville-Caputo) is done too.

## I. INTRODUCTION

The definition of ordinary $s$-order derivative $D^s$ and integral $J^s$ implies that $s$ is a positive integer. May we extend the definition of $D^s$ and $J^s$ for any value of $s$? ($s \in \mathbb{R}, s \in \mathbb{C}$) This is the problem of fractional derivatives first raised by Leibnitz in 1695.:"Can the meaning of derivatives with integer order by generalized to derivatives with non-integer orders?" L'Hospital replied to Liebnitz by another question:"What if the order is 1/2?" And Leibnitz responded: "It will lead to a paradox, from which one day useful consequences will be drawn". This was the birth of fractional derivatives. Since then several approaches have been done [1], [2],[3].

The strategy to tackle the problem is the following :
1) define $J^1(f)(x)$ and $D^1(f)(x)$
2) define from $J^1(f)(x)$ and $D^1(f)(x)$   $J^s(f)(x)$ and $D^s(f)(x)$ for all $s \in \mathbb{N}$
3) extend to $s \in \mathbb{N}$, then to $s \in \mathbb{R}$
4) extend to $s \in \mathbb{C}$
5) and finally, look for in which case the results satisfy the following relations.

5a) Principle of correspondance

$$\lim_{s \to n} D^s(f)(x) = \frac{d^n}{dx^n} f(x) \quad for\ n \in \mathbb{N}$$

$$\lim_{s \to n} J^s(f)(x) = \int_a^x \int_a^{t_1} \int_a^{t_2} \ldots \int_a^{t_n} f(t_n) dt_n dt_{n-1} \ldots dt_1$$

$for\ n \in \mathbb{N}$

5b) linearity property of $D^s$ and $J^s$

$$D^s(cf)(x) = cD^s(f)(x)$$
$$D^s(f_1 + f_2)(x) = D^s(f_1)(x) + D^s(f_2)(x)$$

or equivalently

$$D^s(c^1 f_1 + c^2 f_2)(x) = c^1 D^s(f_1)(x) + c^2 D^s(f_2)(x)$$

for any constants $c^1$ et $c^2$ and for any causal functions $f_1$ and $f_2$ (the same relations for $J^s(f)(x)$)

5c) Is the value of the fractional derivative $D^s$ of a constant function null or not?
5d) Semi-group property

$$J^{s_1}(J^{s_2})(f)(x) = J^{s_1+s_2}(f)(x) = J^{s_2}(J^{s_1})(f)(x)$$
$$D^{s_1}(D^{s_2})(f)(x) = D^{s_1+s_2}(f)(x) = D^{s_2}(D^{s_1})(f)(x)$$

or

$$J^{s_1} J^{s_2} = J^{s_1+s_2} = J^{s_2} J^{s_1}$$
$$D^{s_1} D^{s_2} = D^{s_1+s_2} = D^{s_2} D^{s_1}$$

This is the strategy that we are going to apply. As extensively admitted, we use the adjective fractional in a broad meaning: fractional, real, complex and the operators product notation

$$ABC(f) = A[B[C(f)]]]$$

## II. DEFINITION OF ONE ORDER INTEGRAL OPERATOR $J^1$ AND ONE ORDER DERIVATIVE $D^1$

Let $E$ be the set of integrable and derivable causal function $f$ defined on the interval $I = ]a, +\infty[, a \in \mathbb{R}$ such

$$f(x) = 0\ for\ x \leq a$$

$f(a)$ may be zero ($f(a) = 0$) or different of zero.
We define the one order integral operator $J^1$ and the one order derivative operator $D^1$ by the relations

$$J^1(f)(x) = \int_a^x f(t)dt \qquad (2.1)$$

$$D^1(f)(x) = \frac{d}{dx}(f)(x) = f^{(1)}(x) \qquad (2.2)$$

$f^{(1)}$ being the first derivative of $f$.

## III. RELATION BETWEEN $J^1$ AND $D^1$ OPERATORS

The following relations may be shown easily

$$D^1 J^1(f)(x) = D^1 \int_a^x f(t)dt = D^1[F(x) - F(a)] = f(x) \quad (3.1)$$

$F(x)$ is a primitive of $f(x)$.

$D^1 J^1 = 1_E$ where $1_E$ is the identity operator over $E$.
If $f(a) = 0$,

$$J^1 D^1(f)(x) = J^1(D^1 f)(x) = \int_a^x f^{(1)}(t)dt$$
$$= f(x) - f(a) = f(x)$$

$$J^1 D^1 = D^1 J^1 = 1_E \quad (3.3a)$$

If $f(a) \neq 0$

$$J^1 D^1 \neq 1_E \quad (3.3b)$$

From the theory of the inverse of an operator, any operator $A$ has *at least* a right handside inverse $A_R$ such as

$$AA_R = 1_{Val(A)} \quad (3.4)$$

where $Val(A)$ is the value domain of $A$ and $A_R$ depends on a choice [4], [5]
By an appropriate mechanism, $A_R$ may be chosen to be unique in order that $A_R$ is an operator.
Taking account of the relation (3.3a) and (3.3b), we may define for $A \equiv J^1$ the right handside inverse $D_R^1$ of $J^1$

$$J^1 D_R^1 = 1_{Val(J^1)} \quad (3.5)$$

with the choice of $f(a) = 0$.

**Remarks**

a) The condition $f(a) = 0$ depends on two parameters : the function $f$ and the choice of the lower bound $a$ in the integral $J^1$ (2.1)
For instance,
if $f(x) = e^x$, $a$ may be taken equal to $-\infty$,
if $f(x) = e^{-x}$, $a$ may be taken equal to $+\infty$,
if $f(x) = P_k(x)$ is a polynom of degree k, the lower bound $a$ may be choosen equal to zero and $f(0) = 0$ implies that the term independent on $x$ is equal to zero too.
If $f(x)$ is equal to a constant $C$, the condition $f(0) = 0$ implies that the constant $C = 0$.

b) We assume $f(a) = 0$ then $D_R^1 = J^{-1}$ and we have
$$D^1 J^1 = J^1 D_R^1 = 1_E \quad (3.6)$$

## IV. DEFINITIONS OF $J^s(f)(x)$ FOR ANY VALUE OF THE NUMBER $s$

For any positive integer $s$, we define $J^s$ by iterating $s$-times $J^1(f)(x)$

$$J^0(f)(x) = f(x) \text{ or } J^0 = 1_E$$

$$J^1(f)(x) = \int_a^x f(t)dt$$

$$J^2(f)(x) = \int_a^x J^1 f(t)dt = \int_a^x [\int_a^{t_1} f(t_2) \, dt_2] \, dt_1$$

...

$$J^s(f)(x) = \int_a^x J^{s-1} f(t)dt$$
$$= \int_a^x \int_a^{t_1} \int_a^{t_2} ... \int_a^{t_{s-1}} f(t_s) \, dt_s dt_{s-1} ... dt_2 dt_1 \quad (4.1)$$

It may be shown (see for instance [5])

$$J^s(f)(x) = \frac{1}{\Gamma(s)} \int_0^x (x-y)^{s-1} f(y) dy \quad (4.2)$$

$$= \frac{x^s}{\Gamma(s)} \int_{a/x}^1 (1-u)^{s-1} f(ux) du \quad (4.3)$$

where $\Gamma(s)$ is the Euler's gamma function for $s$. We have extended the definition of $J^s$ for any $s \in \mathbb{R}_+$ [5] and for any complex number $s \in \mathbb{C}$ with the condition $Re(s) > 0$ [6] and studied their properties in the case of $f(a) = 0$.

## V. SOME PROPERTIES OF $J^s$

*Theorem 1*
$J^s$ is a linear operator over $E$.

*Proof*

$$J^s(\lambda f + \mu g)(x) = \lambda J^s(f)(x) + \mu J^s(g)(x) \quad (5.1)$$

for $f$ and $g$ belonging to $E$ and for any complex $\lambda$ and $\mu$. The proof is very easy.

*Theorem 2: Semi-group property of $J^s$*

$$J^{s_1} J^{s_2} = J^{s_1+s_2} = J^{s_2} J^{s_1} \quad (5.2)$$

*Proof*
We have to show

$$J^{s_1} J^{s_2}(f)(x) = J^{s_1+s_2}(f)(x)$$

for any $f \in E$ and for any $s_1$ and $s_2$, real integers, real fractional numbers, complex numbers with $Re(s_1) > 0$ and $Re(s_2) > 0$.
Let us suppose first $s_1$ and $s_2$ real integer numbers

$$J^{s_1}J^{s_2}(f)(x) = \frac{1}{\Gamma(s_1)\Gamma(s_2)} \int_a^x (x-y)^{s_1-1} dy \int_a^y (y-z)^{s_2-1} f(z) dz$$

We apply the Dirichlet's formula given by Whittaker and Watson [7],[8]

$$\int_a^x dy(x-y)^{\alpha-1} \int_a^y dz\, (y-z)^{\beta-1} dy g(y,z)$$

$$= \int_a^z dz \int_z^x dy(x-y)^{\alpha-1}(y-z)^{\beta-1} g(y,z) \quad (5.3)$$

for $\alpha = s_1, \beta = s_2, g(y,z) = f(z)$

$$J^{s_1}J^{s_2}(f)(x) = \frac{1}{\Gamma(s_1)\Gamma(s_2)} \int_a^x dz\, f(z) \int_z^y dy(x-y)^{s_1-1}(y-z)^{s_2-1}$$

Then, we change the variable $y$ in the second integral into $u$

$$u = \frac{y-z}{x-z}$$

or
$$y = z + u(x-z)$$

$$dy = du(x-z)$$

$$\int_z^y dy(x-y)^{s_1-1}(y-z)^{s_2-1}$$

$$= (x-z)^{s_1+s_2-1} \int_0^1 du(1-u)^{s_1-1} u^{s_2-1}$$
$$= (x-z)^{s_1+s_2-1} B(s_1)(s_2)$$

where $B$ is the Euler's beta function

$$B(s_1)(s_2) = \frac{\Gamma(s_1)\Gamma(s_2)}{\Gamma(s_1+s_2)}$$

then

$$J^{s_1}J^{s_2}(f)(x) = \frac{1}{\Gamma(s_1+s_2)} \int_a^x dz f(z)(x-z)^{s_1+s_2-1}$$

$$= J^{s_1+s_2}(f)(x)$$

As $J^{s_1+s_2}$ is symmetric in $s_1$ and $s_2$, we have

$$J^{s_1}J^{s_2}(f)(x) = J^{s_1+s_2}(f)(x) = J^{s_2}J^{s_1}(f)(x)$$

or

$$J^{s_1}J^{s_2} = J^{s_1+s_2} = J^{s_2}J^{s_1}$$

We have supposed $s_1$ and $s_2$ are positive integer numbers. But if we use the extensions of Euler's gamma and beta functions, it may be easily shown that the semi-group property of $J^s$ stands true for any $s_1$ and $s_2$ (positive fractional, positive real,) [5] and for any $s_1$ and $s_2$ complex numbers with the conditions $Re(s_1) > 0$ and $Re(s_2) > 0$ [6].

This important property is very useful to shorten the demonstration of many interesting formulae, in particular the semi-group property of $s$-order derivative $D^s$. (see section VI)

## VI. DEFINITION OF THE $s$-ORDER DERIVATIVE OPERATOR $D^s$ FROM $J^s$

We deduce the definition of the $s$−order derivative operator $D^s$ from the definitions of $s$-order integral operator $J^s$.
For this purpose, we take advantage of the relation (3.3) or (3.6)

$$J^1 D^1 = D^1 J^1 = 1_E$$

$D^1$ is the inverse operator $J^{-1}$ of $J^1$  (6.1)
$J^1$ is the inverse operator $D^{-1}$ of $D^1$  (6.2)

By mathematical induction, for $s$ positive integer, we may define

$$D^s = J^{-s} \quad (6.3)$$

$$J^s = D^{-s} \quad (6.4)$$

Then we may extend the definition of $J^s$ (see formulae 4.2 and 4.3) given for $s$ positive integer into $s$ negative integer for any real (positive or negative number) and for any complex number by using the semi-group property of $J^s$.
For the definition of $D^s$ for $s > 0$, we have two possibilities left-hand derivative $D_L^s$ and right-hand derivative $D_R^s$:

a) $D_L^s = J^{-k}J^{k-s} = D^k J^{k-s}$ for $k \in \mathbb{N}$  (6.5)

$$D_L^s(f)(x) = D^k J^{k-s}(f)(x) \quad (6.6)$$
$$= D^k \frac{1}{\Gamma(k-s)} \int_a^x (x-y)^{k-s-1} f(y) dy \quad (6.7)$$
$$= D^1 \frac{1}{\Gamma(1-s)} \int_a^x (x-y)^{-s} f(y) dy \quad (6.8)$$

with the choice $k = 1$ in the relation (6.7).
If $f(x) = C$ where $C$ is a constant.

$$D_L^s(C) = \frac{C(x-a)^{-s}}{\Gamma(1-s)} \neq 0 \quad C \neq 0 \quad (6.9)$$

b) $D_R^s = J^{k-s}J^{-k} = J^{k-s}D^k$ for $k \in \mathbb{N}, k > s$  (6.10)

$$D_R^s(f)(x) = J^{k-s}(D^k f)(x) \quad (6.11)$$

$$= \frac{1}{\Gamma(k-s)} \int_a^x (x-y)^{k-s-1} f^{(k)}(y) dy \quad (6.12)$$

$$= \frac{1}{\Gamma(1-s)} \int_a^x (x-y)^{-s} f^{(1)}(y) dy \quad (6.13)$$

with the choice $k = 1$.
If $f(x) = C$ where $C$ is a constant,

$$D_R^s C = 0 \text{ even if the constant } C \neq 0 \quad (6.14)$$

The choice $D_R^s$ in the relation (6.10) is better than $D_L^s$ in the relation (6.5) because of the relation (6.14) instead of the relation (6.9) which does not verify the condition: the fractional derivative of a constant is equal to zero.

c) If $0 < s < 1$, for simplicity sake, we may choose $k = 1$

$$D^s(f)(x) = D^1 \frac{1}{\Gamma(1-s)} \int_a^x (x-y)^{-s} f(y) dy \quad (6.15)$$

$$= \frac{1}{\Gamma(1-s)} \frac{d^1}{dx^1} \int_a^x (x-y)^{-s} f(y) dy \quad (6.16)$$

We apply the derivative under the integration sign and we integrate by part after, the terms $(x-y)^{-s} f(y)$ to be taken for $y = s$ minus one term for $y = a$ cancels each other without introducing any conditions on the function $f$. Finally we obtain

$$D^s(f)(x) = \frac{1}{\Gamma(1-s)} \int_a^x (x-y)^{-s} f^{(1)}(y) dy \quad (6.17)$$

$$= \frac{1}{\Gamma(1-s)} \int_a^x (x-y)^{-s} f^{(1)}(y) dy \quad (6.18)$$

The two expressions (6.17) and (6.18) are exactly the Liouville-Caputo fractional derivative definition [7] if we take $a = -\infty$:

$$_{LC}D_+^s(f)(x) = {}_L J_+^{1-s} \frac{d}{dx}(f)(x) \quad (6.19)$$

It consists to change $s$ into $-s$ in the expression (4.2) of $J^s$. The lower script LC means Liouville-Caputo; the lower scripts L and + in $_L J_+^{1-s}$ mean respectively Liouville for Liouville fractional integral definition, the integral being "left handed integral"; it collects weighted function values for $y < x$ which means left hand side from $x$. If $y$ is time-like coordinate, the left handed integral corresponds to a "causal" function $f$. In our notations, the relation (6.15) for $k = 1$ is

$$D^{(s)}(f)(x) = J^{1-s} D^1(f)(x) \quad (6.20)$$

$$= J^{1-s} D^1(f)(x) \quad (6.21)$$

d) If we take $k = 2, 3, 4, \ldots$ we will obtain the same result; direct calculations are not difficult but tedious.

e) If $s > 1$, we put $s = E(s) + s_1$ where $E(s)$ is the entire part of $s$ and $0 \leq s_1 < 1$.

$$D^s = D^{E(s) + s_1}$$

$$= D^{s_1} D^{E(s)} \; (semi-group \; property \; of \; D^s)$$

It means that the fractional derivative is an ordinary standard derivative $D^{E(s)}$ followed by a fractional derivative $D^{s_1}$ with $0 \leq s_1 < 1$.

## VII. DUAL PROPERTIES OF FRACTIONAL DERIVATIVES

It is easy to derive from properties of fractional integrals $J^s$ the following properties for $D^s$

- linear property

$$D^s(\lambda f + \mu g)(x) = \lambda D^s(f)(x) + \mu D^s(g)(x) \quad (7.1)$$

for $\lambda, \mu \in \mathbb{C}$ and $f, g \in E$

- semi-group property

$$D^{s_1} D^{s_2} = D^{s_1 + s_2}$$
$$= D^{s_2} D^{s_1} \quad (7.2)$$

The direct proofs are not difficult, they may be performed by a similar way as for fractional integral calculations

## VIII. LIMITS FOR FRACTIONAL INTEGRALS AND DERIVATIVES

Let us now look for the correspondence principle. We have to show the following limits

$$\lim_{s \to n} J^s(f)(x) = J^n(f)(x) \quad (s \in \mathbb{R}_+, n \in \mathbb{N})$$

$$\lim_{s \to n} D^s(f)(x) = D^n(f)(x) \quad (s \in \mathbb{R}_+, n \in \mathbb{N})$$

where $J^n$ and $D^n$ are respectively ordinary integrals and ordinary derivatives.

*Theorem 3*

Let us assume first $0 \leq s \leq 1$

$$\lim_{s \to 0^+} J^s(f)(x) = f(x) \; \forall f \in E$$

$$\lim_{s \to 0^+} D^s(f)(x) = f(x) \; \forall f \in E$$

or

$$\lim_{s \to 0^+} J^s = J^0 = 1_E$$

$$\lim_{s \to 0^+} D^s = D^0 = 1_E$$

*Proof*

$$J^s(f)(x) = \frac{1}{\Gamma(s)} \int_a^x (x-y)^{s-1} f(y) dy \quad (8.1)$$

$$= \frac{1}{\Gamma(s+1)} \int_a^x (x-y)^s f^{(1)}(y) dy \quad (8.2)$$

where

$$f^{(1)}(y) = \frac{d}{dy} f(y)$$

For $s = \varepsilon \ll 1$

$$J^\varepsilon(f)(x) = \frac{1}{\Gamma(1+\varepsilon)} \int_a^x (x-y)^\varepsilon f^{(1)}(y) dy$$

We give in the appendix the Taylor's series expansion of

$$\frac{1}{\Gamma(1+\varepsilon)}(x-y)^\varepsilon = 1 + \varepsilon[\gamma + ln(x-y)]$$
$$+ \frac{\varepsilon^2}{2}[\gamma^2 - \frac{\pi^2}{6} + 2\gamma ln(x-y) + ln^2(x-y)] + \cdots \quad (8.3)$$

where $\gamma = 0.577\,215\,665\,901\,532\,860\,6\ldots$ is the Euler – Mascheroni's constant [7]

$$J^0(f)(x) = \int_a^x f^{(1)}(y) dy \quad (8.4)$$
$$= f(x) \quad (8.5)$$

because $f(a) = 0$
or $J^0 = 1_E$
By duality, we have the same relations for $D^s$.

*Theorem 4*

We have the dual relations for $J^s$ and $D^s$, $s \in \mathbb{C}$ with $Re(s) > 0$.

*Proof*
The proof is obtained easily in copying the demonstration in the real case.

*Theorem 5*

For instance, let us verify directly the relations

$$\lim_{\varepsilon \to 0^+} D^{1-\varepsilon} f(x) = \frac{d}{dx} f(x) \quad (8.7)$$

*Proof*
We start from the expression (6.15) of $D^s(f)(x)$ for $s = 1 - \varepsilon$ with $\varepsilon > 0$ and $\varepsilon \ll 1$

$$D^{1-\varepsilon}(f)(x) = \frac{1}{\Gamma(\varepsilon)} \int_a^x (x-y)^{\varepsilon-1} \frac{df}{dx} dx \quad (8.8)$$

We integrate first by parts

$$D^{1-\varepsilon}(f)(x) = \frac{1}{\Gamma(\varepsilon)}\{-[\frac{(x-y)^\varepsilon}{\varepsilon} f^{(1)}(y)]_{y=a}^{y=x}$$
$$+ \int_a^x \frac{(x-y)^\varepsilon}{\varepsilon} f^{(2)}(y) dy\}$$

We integrate by parts again the integral containing $f^{(2)}$. The integrated part cancels exactly the first term in the right hand side of the relation

$$D^{1-\varepsilon}(f)(x) = \frac{1}{\Gamma(\varepsilon)} \int_a^x \frac{(x-y)^\varepsilon}{\varepsilon} f^{(2)}(y) dy \quad (8.10)$$
$$= \frac{1}{\Gamma(1+\varepsilon)} \int_a^x (x-y)^\varepsilon f^{(2)}(y) dy \quad (8.11)$$

We utilise the Taylor's series expansion (8.3)

$$\lim_{\varepsilon \to 0^+} D^{1-\varepsilon}(f)(x) = \int_a^x f^{(2)}(y) dy$$
$$= f^{(1)}(x) - f^{(1)}(a)$$
$$= \frac{d}{dx} f(x) \quad \text{if and only if } f^{(1)}(a) = 0$$

*Theorem 6*

$$\lim_{\varepsilon \to 0^+} J^\varepsilon(f)(x) = \int_a^x f^{(1)}(y) dy$$
$$= f(x) - f(a)$$
$$= f(x) \quad \text{if and only if } f(a) = 0$$

$$\lim_{\varepsilon \to 0^+} J^\varepsilon = 1_E \qquad \lim_{\varepsilon \to 0^-} J^\varepsilon = 1_E$$

*Theorem 7*

$$\lim_{\varepsilon \to 0^+} J^{1-\varepsilon}(f)(x) = J^1 \quad (8.12)$$

*Proof*
$$\lim_{\varepsilon \to 0^+} J^{1-\varepsilon}(f)(x) = J^1 \lim_{\varepsilon \to 0^+} J^{-\varepsilon}(f)(x)$$
$$= J^1(f)(x)$$

because

$$\lim_{\varepsilon \to 0^+} J^{-\varepsilon}(f)(x) = f(x)$$

The direct calculation is not difficult.

$$J^{1-\varepsilon}(f)(x) = \frac{1}{\Gamma(1-\varepsilon)} \int_a^x (x-y)^{-\varepsilon} f(y) dy$$

We utilise the Taylor's series expansion (8.3) in changing $\varepsilon$ into $-\varepsilon$

$$\lim_{\varepsilon \to 0^+} J^{1-\varepsilon}(f)(x) = \int_a^x f(y) dy = J^1(f)(x)$$

$$\lim_{\varepsilon \to 0^+} J^{1-\varepsilon} = J^1$$

*Theorem 8.*

Let us assume now $s > 1$ with $n = E(s) \in \mathbb{N}$
We put $s = n + s_1$ with $0 \leq s_1 < 1$

$$\lim_{s_1 \to 0^+} J^{n+s_1}(f)(x) = J^n(f)(x) \quad (8.13)$$
$$\lim_{s_1 \to 0^+} J^{n+1-s_1}(f)(x) = J^{n+1}(f)(x) \quad (8.14)$$

Proof
The proof is immediate by utilizing the semi-group property of $J^s$.

**Remark**

The relation (8.7) may be obtained by duality from the relation (5.2) by taking account of

$$D^{1-\varepsilon}(f)(x) = J^{\varepsilon-1+1}(f^{(1)})(x)$$

$$\lim_{\varepsilon \to 0^+} D^{1-\varepsilon} = \lim_{\varepsilon \to 0^+} J^{\varepsilon}(f^{(1)})(x)$$
$$= f^{(1)}(x) = \frac{d}{dx}f(x)$$

## IX. CASE OF $f^{(m)}(a) \neq 0$ for $m = 0, 1, 2, 3, ...$

In this paragraph, let us assume $f^{(m)}(a) \neq 0$ for $m = 0, 1, 2, 3$ ... and look for some consequences.

*Theorem 9*

$$J^s(f)(x) = J^{s+1}(f^{(1)})(x) + \frac{1}{\Gamma(s+1)}(x-a)^s f^{(0)}(a) \quad (9.1)$$

where

$$f^{(1)} = \frac{d}{dx}(f) = D^1(f)$$

is the first ordinary derivative of $f$

*Proof*

$$J^s(f)(x) = \frac{1}{\Gamma(s)}\int_a^x (x-y)^{s-1} f(y) dy$$

We integrate by parts and take account of $s\Gamma(s) = \Gamma(s+1)$

$$J^s(f)(x) = J^{s+1}(f^{(1)})(x) + \frac{(x-a)^s}{\Gamma(s+1)}f(a)$$

*Theorem 10.*

Theorem 9 may be generalized for any positive integer $k$

$$J^s(f)(x) = J^{s+k}(f^{(k)})(x) + \sum_{m=1}^{m \leq k} \frac{(x-a)^{s+m-1}}{\Gamma(s+m)} f^{(m-1)}(a) \quad (9.2)$$

$$= \frac{1}{\Gamma(s+k)}\int_a^x (x-y)^{s+k-1} f^{(k)}(y) dy$$

$$+ \sum_{m=1}^{m \leq k} \frac{(x-a)^{s+m-1}}{\Gamma(s+m)} f^{(m-1)}(a)$$

*Proof*

We utilize mathematical induction. The relation is true for $k = 1$. Let us suppose that it is true for $k$ and we will show that it stands true for $k + 1$.

$$J^s(f)(x) = \frac{1}{\Gamma(s+k)}\int_a^x (x-y)^{s+k-1} f^{(k)}(y) dy$$

$$+ \sum_{m=1}^{m \leq k} \frac{(x-a)^{s+m-1}}{\Gamma(s+m)} f^{(m-1)}(a)$$

We integrate by parts the integral in the right hand side

$$\frac{1}{\Gamma(s+k)}\int_a^x (x-y)^{s+k-1} f^{(k)}(y) dy =$$
$$\frac{1}{\Gamma(s+k)}\int_a^x \frac{(x-y)^{s+k}}{(s+k)} f^{(k+1)}(y) dy$$
$$+ \frac{1}{\Gamma(s+k)}\frac{(x-a)^{s+k}}{(s+k)} f^{(k)}(a)$$

Then

$$J^s(f)(x) = \frac{1}{\Gamma(s+k+1)}\int_a^x (x-y)^{s+k} f^{(k+1)}(y) dy$$

$$+ \sum_{m=1}^{m \leq k+1} \frac{(x-a)^{s+m-1}}{\Gamma(s+m)} f^{(m-1)}(a)$$

## X. COMPARISON WITH OTHER DEFINITIONS

a) If $f^{(m)}(a) = 0$ for any $m = 0, 1, 2, 3, ...$ then Theorem 7 and Theorem 8 give

$$J^s(f)(x) = J^{s+1}(f^{(1)})(x) \quad (10.1)$$

$$= J^{s+k}(f^{(k)})(x) \quad (10.2)$$

Or

$$J^s(f)(x) = J^{s+k}(D^k f)(x) \quad (without\ summation\ over\ k)$$

$$J^s = J^{s+k} D^k \quad (10.3)$$

These relations may be derived from the semi-group property of $J^s$ (5.1) and $D^k = J^{-k}$

b) If $a = -\infty$, the relation (4.2) gives

$$J^s(f)(x) = \frac{1}{\Gamma(s)}\int_{-\infty}^x (x-y)^{s-1} f(y) dy \quad (10.4)$$

which is exactly the definition of Liouville fractional integral $_L J_+^s(f)(x)$ ; the underscript $L$ stands for Liouville and the underscript $(+)$ is related to the left-handed integral [formula (5.15), p.36 in reference [7].

c) If $a = 0$, the relation (4.2) gives

$$J^s(f)(x) = \frac{1}{\Gamma(s)}\int_0^x (x-y)^{s-1} f(y) dy \quad (10.5)$$

which is the Riemann fractional integral $_R J_+^s(f)(x)$ for left handed integral (formula (5.17) p.36 in reference [7])

d) For the expression of the fractional derivative $D^s$, we have to change the sign of $s$ in the expression of the fractional integral $J^s$ in the relation (4.2) with $k = 1$.

$$D^s(f)(x) = D^1 J^{1-s}(f)(x)$$

$$= \frac{d}{dx}\frac{1}{\Gamma(1-s)}\int_a^x (x-y)^{-s} f(y) dy \quad (10.6)$$

If $a = -\infty$, we have

$$D^s(f)(x) = \frac{d}{dx} \frac{1}{\Gamma(1-s)} \int_{-\infty}^{x} (x-y)^{-s} f(y) dy \quad (10.7)$$

which is the definition of Liouville fractional derivative $_L D_+^s(f)(x)$ for a left handed integral for a causal function $f$. [formula (5.29) p.38 of reference [7]]

e) If $a = 0$, the relation (9.8) gives

$$D_R^s(f)(x) = \frac{d}{dx} \frac{1}{\Gamma(1-s)} \int_0^x (x-y)^{-s} f(y) dy \quad (10.8)$$

which is the definition of Riemann fractional derivative $_R D_+^s(f)(x)$ (formula 5.33 p.39 of reference 7)

f) Let us change $s$ into $-s$ in the expression of $J^s(f)(x)$ in the relation (10.1) $(D^s = J^{-s})$:

$$\begin{aligned} D_R^s(f)(x) &= J^{-s}(f)(x) \\ &= J^{1-s}(f^{(1)})(x) \\ &= \frac{1}{\Gamma(1-s)} \int_a^x (x-y)^{-s} \frac{d}{dy} f(y) dy \end{aligned} \quad (10.9)$$

g) If $a = 0$, the relation (10.9) is the definition of Caputo fractional derivative $_C D_+^s(f)(x)$ for a causal function $f$ (formula 5.44 p.42 of the reference [7])

h) If $a = -\infty$, the relation (10.9) is the definition of Liouville-Caputo fractional derivative $_{LC} D_+^s(f)(x)(x)$ for a causal function $f$.

## XI. CONCLUSIONS

Let us conclude by stressing once more that the definition of fractional derivative $D^s$ may be obtained from the definition of an integer order (ordinary) derivative $D^k$ and fractional integral by

$$D^s(f)(x) = D^k J^{k-s}(f)(x)$$

$$D^s = D^k J^{k-s} \text{ for any integer } k > s$$

In the second members of these relations, we point out that there is no summation over $k$ and it is independent on $k$.

Comparing the expression (4.2) of $J^s(f)(x)$ and (6.3), (6.4), (6.5), (6.17) of $D^s(f)(x)$, we obtain an **unified** definition for fractional integrals and fractional derivatives according to **the sign of the real part** of the order $s$.

## APPENDIX

Let us calculate the Taylor's series expansion (8.3)

$$\phi(\varepsilon) = \frac{1}{\Gamma(1+\varepsilon)} (x-y)^\varepsilon$$

$$= \frac{1}{\Gamma(1+\varepsilon)} e^{\varepsilon \ln(x-y)}$$

$$= \phi(0) + \frac{\varepsilon}{1!} \left(\frac{\partial \phi}{\partial \varepsilon}\right)_{\varepsilon=0} + \frac{\varepsilon^2}{2!} \left(\frac{\partial \phi}{\partial \varepsilon}\right)_{\varepsilon=0} + \cdots$$

The calculation may be simplified by utilizing Leibnitz's formula

$$\left(\frac{a^\varepsilon}{\Gamma(1+\varepsilon)}\right)^{(n)} = \sum_{k=1}^{n} \binom{k}{n} \left(\frac{1}{\Gamma(1+\varepsilon)}\right)^{(k)} (a^\varepsilon)^{(n-k)}$$

with

$$a^\varepsilon = \exp(\varepsilon \ln a) \quad a = x - y$$

For $n = 0$

$$\phi(0) = \frac{1}{\Gamma(1)} a^0 = 1$$

For $n = 1$

$$\left(\frac{a^\varepsilon}{\Gamma(1+\varepsilon)}\right)^{(1)} = -\frac{\Gamma^{(1)}(1+\varepsilon)}{[\Gamma(1+\varepsilon)]^2} a^\varepsilon + \frac{1}{\Gamma(1+\varepsilon)} (\ln(a)) a^\varepsilon$$

$$\left(\frac{a^\varepsilon}{\Gamma(1+\varepsilon)}\right)^{(1)}_{\varepsilon=0} = -\Gamma^{(1)}(1) + \ln(a)$$

$$\left(\frac{(x-y)^\varepsilon}{\Gamma(1+\varepsilon)}\right)^{(1)}_{\varepsilon=0} = \gamma + \ln(x-y)$$

Or

$$\Gamma(x) = \int_0^\infty e^{-t} t^{x-1} dt$$

$$\frac{d}{dx} \Gamma(x) = \int_0^\infty e^{-t} t^{x-1} (\ln t) dt$$

$$\frac{d^n}{dx^n} \Gamma(x) = \int_0^\infty e^{-t} t^{x-1} (\ln t)^n dt \quad n \in \mathbb{N}$$

$$\frac{d}{dx} \Gamma(x) \bigg|_{x=1} = \int_0^\infty e^{-t} \ln t \, dt = -\gamma$$

where $\gamma$ is the Euler-Mascheroni's constant

$$\Gamma^{(2)}(1) = \int_0^\infty e^{-t} (\ln t)^2 dt = \gamma^2 + \frac{\pi^2}{6}$$

For $n = 2$

$$\left(\frac{1}{\Gamma(1+\varepsilon)} a^\varepsilon\right)^{(2)} = \left(\frac{1}{\Gamma(1+\varepsilon)}\right)^{(2)} a^\varepsilon + 2\left(\frac{1}{\Gamma(1+\varepsilon)}\right)^{(1)} (a^\varepsilon)^{(1)}$$

$$+ \frac{1}{\Gamma(1+\varepsilon)} (a^\varepsilon)^{(2)}$$

$$\left(\frac{1}{\Gamma(1+\varepsilon)} a^\varepsilon\right)^{(2)}_{\varepsilon=0} = -\Gamma^{(2)}(1) + 2[\Gamma^{(1)}(1)]^2 + 2\Gamma^{(1)}(1) \ln a$$

$$+ (\ln a)^2$$

$$= -\Gamma^{(2)}(1) + 2\gamma^2 + 2\gamma \ln a + (\ln a)^2$$

Or

$$\Gamma^{(2)}(1) = \gamma^2 + \frac{\pi^2}{6}$$

$$\left(\frac{1}{\Gamma(1+\varepsilon)}a^\varepsilon\right)^{(2)}_{\varepsilon=0} = \gamma^2 - \frac{\pi^2}{6} + 2\gamma \ln a + (\ln a)^2$$

Then

$$\frac{1}{\Gamma(1+\varepsilon)}(x-y)^\varepsilon = 1 + \frac{\varepsilon}{1!}(\gamma + \ln(x-y))$$

$$+ \frac{\varepsilon^2}{2!}\{\gamma^2 - \frac{\pi^2}{6} + 2\gamma \ln(x-y) + [\ln(x-y)]^2\} + \cdots$$

We may calculate the coefficients of $\frac{\varepsilon^3}{3!}, \frac{\varepsilon^4}{4!}$ step by step by the method we have just described.